\documentclass[graybox]{svmult}

\usepackage{amsmath,amssymb}
\usepackage{mathptmx}       
\usepackage{helvet}         
\usepackage{courier}        
\usepackage{type1cm}        
\usepackage{makeidx}         
\usepackage{graphicx}        
\usepackage{multicol}        
\usepackage[bottom]{footmisc}
\graphicspath{{pix/}}
\def\e{{\rm e}}

\def\D{\mathcal{D} }
\def\R{\mathbb{R} }
\def\T{\mathbb{T} }

\makeindex             

\begin{document}

\title*{Variance reduction result for a projected adaptive biasing force method}\titlerunning{PABF method}
\author{Houssam AlRachid and Tony leli\`evre}
\institute{Houssam AlRachid \at Universit\'e Paris-Est Cr\'eteil ,  61 Avenue du G\'en\'eral de Gaulle, 94000 Cr\'eteil, France, \email{houssam.alrachid@u-pec.fr}
\and Tony leli\`evre \at \'Ecole des Ponts ParisTech, Universit\'e Paris Est, 6-8 Avenue Blaise Pascal, Cit\'e Descartes, Marne-la-Vall\'ee,  F-77455, France \email{tony.lelievre@enpc.fr}}
%
%
\maketitle

\abstract*{This paper is committed to investigate an extension of the classical  adaptive biasing force method, which is used to compute the free energy related to the Boltzmann-Gibbs measure and a reaction coordinate function. The issue of this technique is that the approximated gradient of the free energy, called biasing force, is not a gradient. The commitment to this field is to project the estimated biasing force on a gradient using the Helmholtz decomposition. The variance of the biasing force is reduced using this technique, which makes the algorithm more efficient than the standard ABF method. We prove exponential convergence to equilibrium of the estimated free energy, with a precise rate of convergence in function of Logarithmic Sobolev inequality constants.}

\abstract{This paper is committed to investigate an extension of the classical  adaptive biasing force method, which is used to compute the free energy related to the Boltzmann-Gibbs measure and a reaction coordinate function. The issue of this technique is that the approximated gradient of the free energy, called biasing force, is not a gradient. The commitment to this field is to project the estimated biasing force on a gradient using the Helmholtz decomposition. The variance of the biasing force is reduced using this technique, which makes the algorithm more efficient than the standard ABF method. We prove exponential convergence to equilibrium of the estimated free energy, with a precise rate of convergence in function of Logarithmic Sobolev inequality constants.}

\section{Introduction}
  Let us consider the Boltzmann-Gibbs measure :
\begin{equation}\label{gib}
\mu(dx)=Z_{\mu}^{-1}\e^{-\beta V(x)}dx,
\end{equation}
where $x\in \D^N$ denotes the position of $N$ particles in $\D\subset \R^n $ (or the $n$-dimensional torus $\T^n$). The potential energy function $V:\D \longrightarrow \R$ associates with the positions of the particles $x\in \D$, $Z_{\mu}$ is the normalization constant and $\beta$ is a constant proportional to the inverse of the temperature.

 The equilibrium probability measure $\mu$ can be sampled through the Overdamped Langevin Dynamics:
\begin{equation}\label{over}
dX_t=-\nabla V(X_t)dt+\sqrt{\frac{2}{\beta}}dW_t,
\end{equation}
where $X_t\in \mathcal{D}^N$ and $W_t$ is a $Nn$-dimensional standard Brownian motion. Under loose assumptions on $V$, the dynamics $(X_t)_{t\geq 0}$ is ergodic with respect to the equilibrium measure $\mu$.
\medskip

Because of the metastability, trajectorial averages converge very slowly to their ergodic limit. To overcome this difficulty, we focus in this paper on the Adaptive Biasing Force (denoted ABF) method (see \cite{dav:01,hen:04}). In order to introduce the ABF method, we require another ingredient: a  reaction coordinate $\xi$ describing the metastable zones of the dynamics associated with the potential energy $V$. For sake of simplicity, take $\xi: (x_1,\ldots,x_n) \in
{\mathbb T^n}  \mapsto (x_1,x_2) \in {\mathbb T}^2$. The associated free energy:
\begin{equation*}\label{free0}
 A(x_1,x_2)=-\beta^{-1}
\ln (Z_{\Sigma_{(x_1,x_2)}})=-\beta^{-1}
\ln \int_{{\mathbb T}^{n-2}} e^{-\beta V(x)}dx_3\ldots dx_n.
\end{equation*}
The idea of the ABF method is that, for a well chosen $\xi$, the dynamics associated with the potential $V-A\circ\xi$ is less metastable than the dynamics associated with $V$.  The so called mean force $\nabla A(z)$, can be obtained as:
\begin{equation*}\label{mean}
\nabla A(x_1,x_2)=\displaystyle Z_{\Sigma_{(x_1,x_2)}}^{-1}\int_{\mathbb{T}^{n-2}}f(x){\rm e}^{-\beta V}dx_3...dx_n=\mathbb{E}_\mu(f(X)|\xi(X)=(x_1,x_2)),
\end{equation*}
where $f=(f_1,f_2)=(\partial_1 V,\partial_2 V)$. At time $t$, the mean force is approximated by $F_t^i(z)=\displaystyle\mathbb{E}_\mu[f_i(X_t) |\xi(X_t)=(x_1,x_2)]$, which also, under appropriate assumptions, converges exponentially fast to $\nabla A$ (see \cite{tony:07,tony:08,tony:10}). Despite the fact that $F_t$ converges to a gradient, there is no reason why $F_t$ would be a gradient at time $t$. In this paper, we propose an alternative method, where we approximate $\nabla A$, at any time $t$, by a gradient denoted $\nabla A_t$.

\section{Projected adaptive biasing force method (PABF)} 
In this section, we present the PABF method, by reconstructing the mean force from the estimated one used in the ABF method.

\medskip

In practice, $A_t$ is obtained from $F_t$ by solving the Poisson problem:
\begin{equation}\label{poi}
{\rm div}(\nabla A_t\psi^{\xi}(t,.))={\rm div}(F_t\psi^{\xi}(t,.))
\mbox{ on } {\mathbb T}^2,
\end{equation}
where $\psi^\xi(t,\cdot)$ denotes the density of the random variables $\xi(X_t)$
which is the Euler equation associated to the minimization problem: \\
 \begin{equation*}\label{min}
A_t=\displaystyle \mathop{\mbox{argmin}}_{g\in H^1(\mathbb{T}^2)/\mathbb{R}}\int_{\mathbb{T}^2} |\nabla g-F_t|^2.
\end{equation*}
Solving \eqref{poi} amounts to computing the so-called Helmholtz-Hodge decomposition of the vector field $F_t$ as:
\begin{equation*}\label{hel}
F_t\psi^{\xi}=\nabla  A_t\psi^{\xi}+R_t \quad \mbox{ on } \mathbb{T}^2,
\end{equation*} 
with ${\rm div}(R_t)=0$. In the following we denote by 
\begin{equation*}\label{projint}
\nabla A_t=\mathcal{P}_{\psi^{\xi}} (F_t),\mbox{ on } {\mathbb T}^2
\end{equation*}
the projection of $F_t$ onto a gradient. We will study the longtime convergence of the following Projected adaptive biasing force (PABF) dynamics:
\begin{equation}\label{pabf}
\left\lbrace
\begin{aligned}
 dX_t&=-\nabla\big( V-A_t\circ\xi\big)(X_t)dt+\sqrt{2\beta^{-1}}dW_t,\\
\nabla A_t&=\mathcal{P}_{\psi^{\xi}} (F_t),\\
 F_t^i(x_1,x_2)&=\displaystyle\mathbb{E}[\partial_i V(X_t)|\xi(X_t)=(x_1,x_2)],\,i=1,2,
\end{aligned}
\right.
\end{equation}	

Using entropy techniques, we study the longtime behavior of the nonlinear Fokker-Planck equation which rules the evolution of the density of $X_t$ solution to (\ref{pabf}). The following theorem shows exponential convergence to equilibrium of $A_t$ to $A$, with a precise rate of convergence in terms of the Logarithmic Sobolev inequality constants of the conditional measures $d\mu_{\Sigma(x_1,x_2)}=Z_{\Sigma_{(x_1,x_2)}}^{-1}{\rm e}^{-\beta V}dx_3...dx_n$. 
\medskip

The assumptions we need to prove the longtime convergence of the biasing force $\nabla A_t$ to the mean force $\nabla A$ are the following:
\begin{description}
  \item[\textbf{H1}] $V\in C^2(\mathbb{T}^n)$, $\exists \gamma >0,\,\forall \,3\leq j\leq n,\forall x\in \mathbb{T}^n,\max(|\partial_1\partial_j V(x)|,|\partial_2\partial_j V(x)|)\leq \gamma .$
  \item[\textbf{H2}] $V$ is such that $\exists\rho>0$, the conditional probability measures $\mu_{\Sigma(x_1,x_2)}$ satisfy a Logarithmic Sobolev inequality with constant $\rho$.
\end{description}
The proof of the following  main theorem is provided in \cite{alr}.
\begin{theorem}\label{theo1}
Let us assume \textbf{H1} and \textbf{H2}. The biasing force $\nabla A_t$ converges to the mean force $\nabla A$ in the following sense:
$$\displaystyle \exists C>0,\exists \lambda>0, \forall\, t\geq 0,\,\int_{\mathbb{T}^2}|\nabla A_t-\nabla  A|^2\psi^{\xi}(t,x_1,x_2)dx_1dx_2\leq\frac{8C\gamma^2}{\rho}\e^{-\lambda t}.$$
\end{theorem}
Since, numerically, we use Monte-Carlo methods to approximate $F_t$ and $\nabla A_t$, the variance is an important quantity to assess the quality of the result. The following second main result is a variance reduction result and proved in \cite{alr}. 
\begin{proposition}\label{vared}
For any time $t>0$, the variance of $\mathcal{P}_{\psi^{\xi}} (F_t)$ is smaller than the variance of $F_t$ in the sense:
$$\displaystyle \forall t>0,\, \int_ {\mathbb{T}^2}{\rm Var}(\mathcal{P}_{\psi^{\xi}} (F_t))\leq \int_ {\mathbb{T}^2}{\rm Var}(F_t),$$
where $\mathrm{Var}(F_t)=\mathbb{E}(|F_t|^2)-\mathbb{E}(|F_t|)^2$ and $|F_t|$ being the Euclidian norm.
\end{proposition}

\section{Numerical experiments}\label{num}
This section is devoted to a numerical illustration of the practical value of the projected ABF compared to the standard ABF approach. 
\medskip

 We consider a system composed of $100$ particles in a two-dimensional periodic box. Among these particles, three particles are designated to form a trimer, while the others are solvent particles. All particles interact through several potential functions such as the Lennard-Jones potential, the double-well potential and a potential on the angle formed by the trimer. We choose the reaction coordinate to be the transition from compact to stretched state in each bond of the trimer. We apply now ABF and PABF dynamics to the trimer problem described above. One can refer to \cite{alr} for more detailed descriptions of the model and the used ABF and PABF algorithms.

 \medskip

\begin{figure}[htp]
  \centering
    \includegraphics[width=100mm]{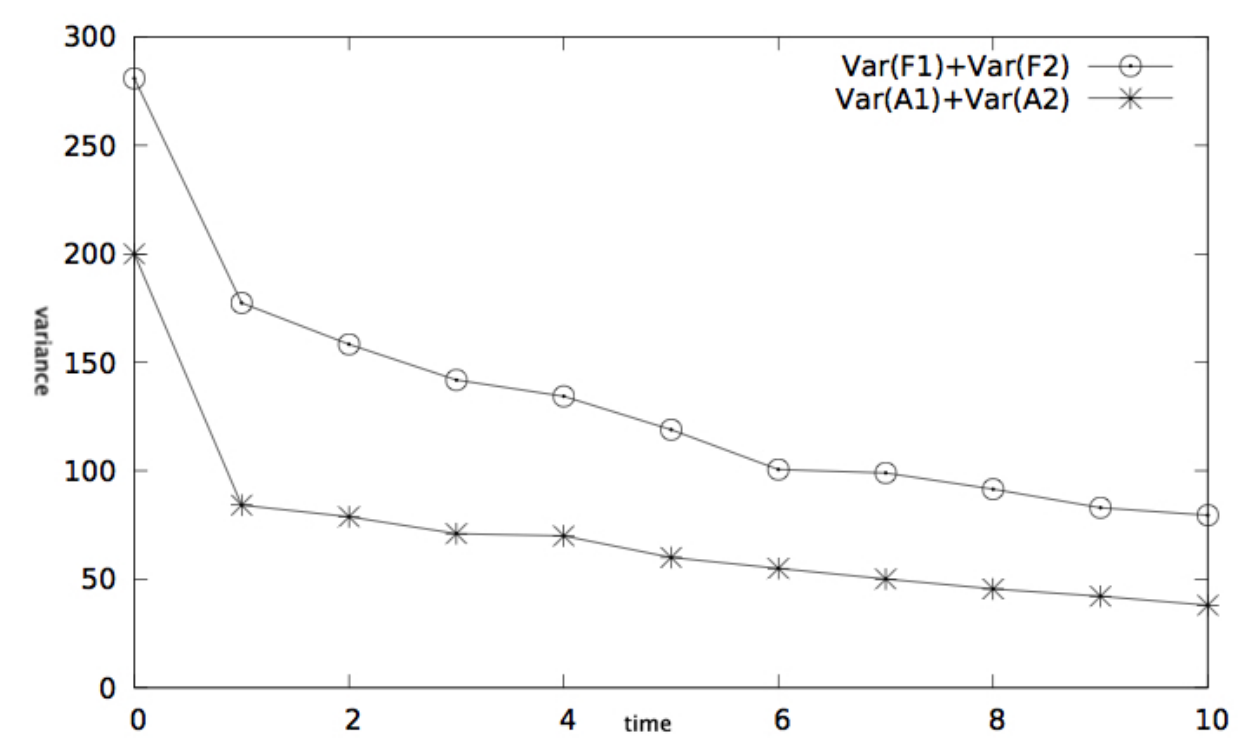}
\caption{Variances as a function of time.}\label{var}
\end{figure}

First, we illustrate the improvement of the projected ABF method in terms of the variances of the biasing forces by comparing $\displaystyle \int{\rm Var}(\nabla A_t)=\int{\rm Var}(\partial_1 A_t)+\int{\rm Var}(\partial_2 A_t)$ (for the PABF method) with $\displaystyle \int{\rm Var}(F_t)=\int{\rm Var}(F_t^1)+\int{\rm Var}(F_t^2)$ (for the ABF method). Figure~\ref{var} shows that the variance for the projected ABF method is smaller than for the standard ABF method, where $ \int{\rm Var}(\nabla A_t)$ (respectively $ \int{\rm Var}(F_t)$) is represented by $Var(F1)+Var(F1)$ (respectively $Var(A1)+Var(A1)$). 

\begin{figure}[htp]
  \sidecaption
    \includegraphics[width=100mm]{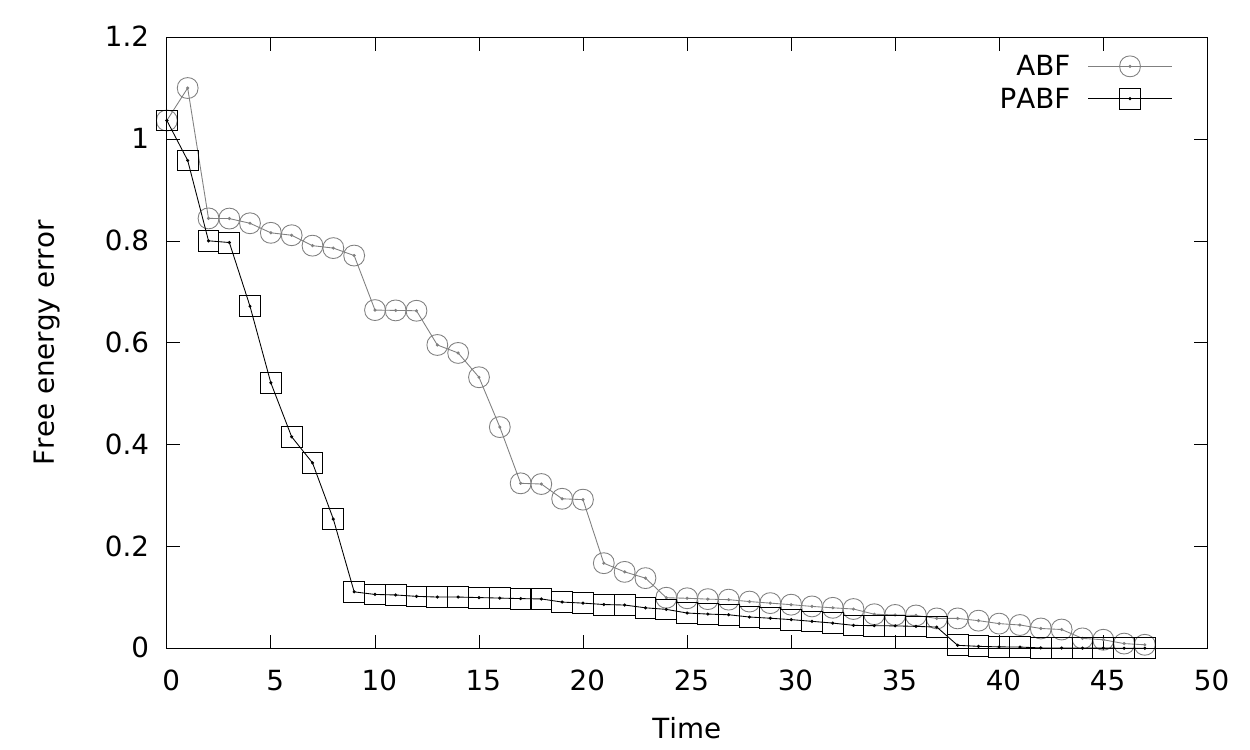}
\caption{ Free energy error as a function of time.}\label{err}
\end{figure}

 We now present, the variation, as a function of time, of the normalized averages $L^2$-distance between the real free energy and the estimated one. As can be seen in Figure~\ref{err}, in both methods, the error decreases as time increases. Moreover, this error is always smaller for the projected ABF method than for the ABF method.

\medskip

\begin{figure}[htp]
  \centering
  \begin{tabular}{cc}
    \includegraphics[width=60mm]{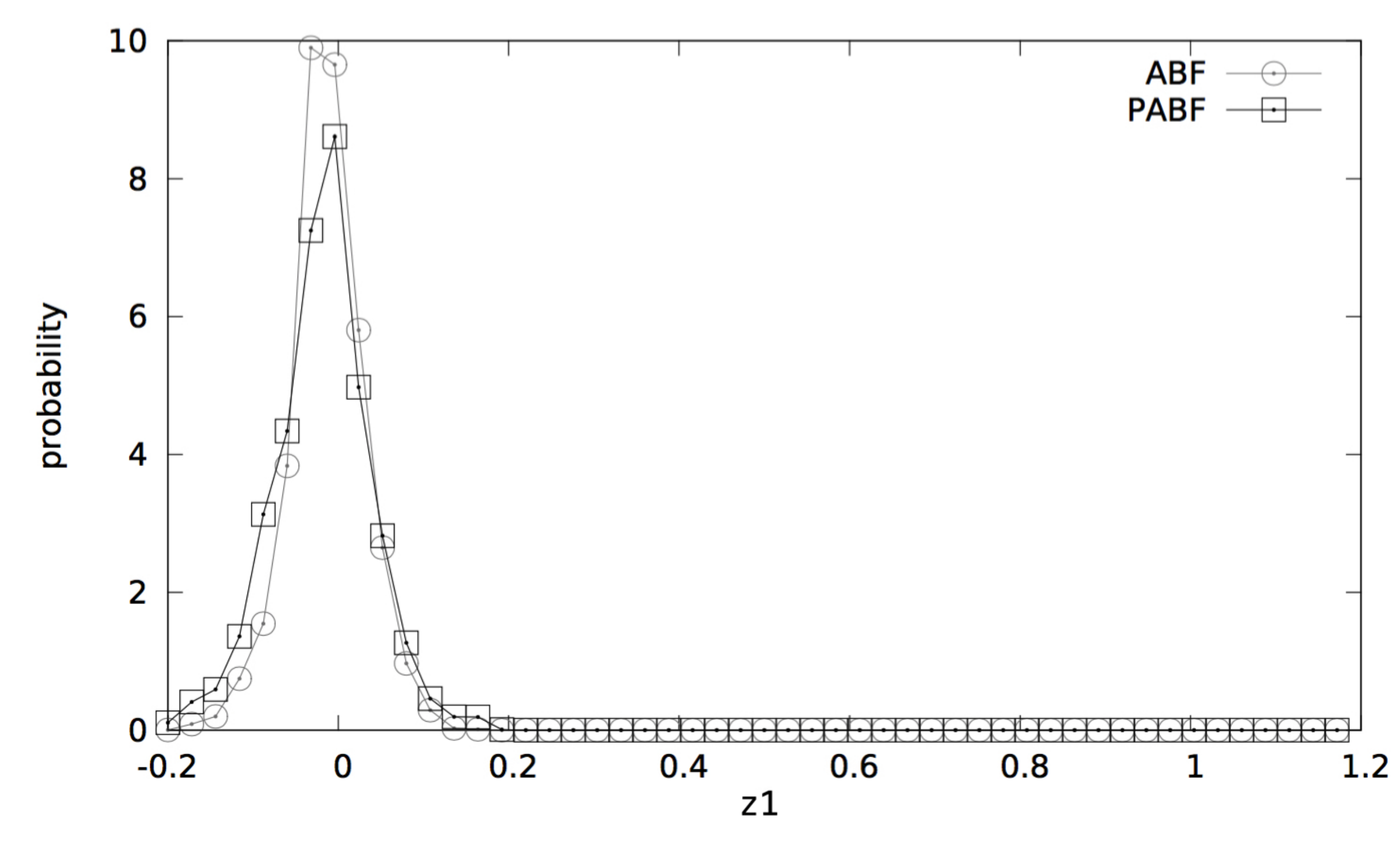}&
    \includegraphics[width=60mm]{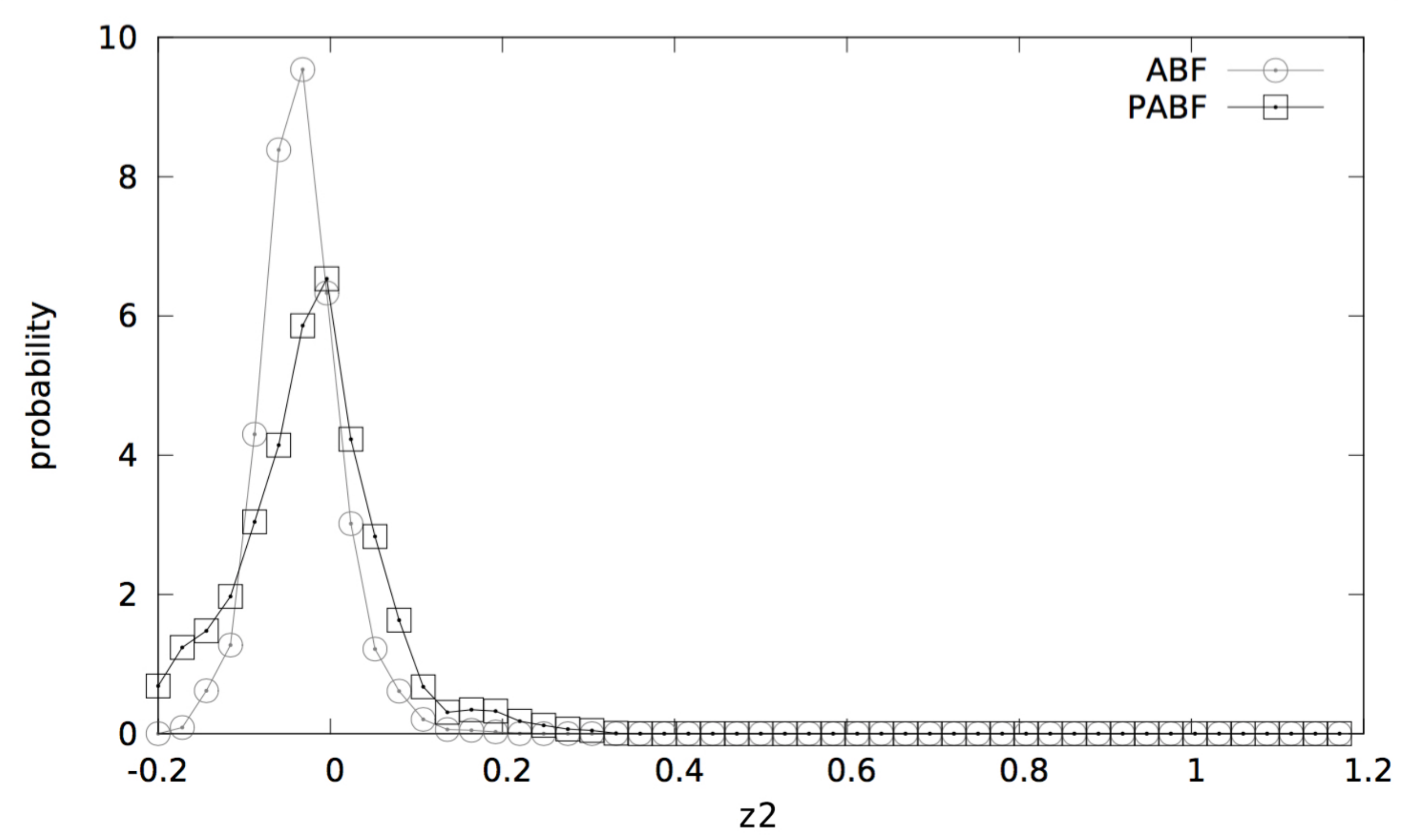}\\
  \end{tabular}
\caption{At time $0.025$. Left: $\int \psi^{\xi}(z_1,z_2)dz_2$; Right:  $\int \psi^{\xi}(z_1,z_2)dz_1$.}\label{dtr0}
\end{figure}

\begin{figure}[htp]
  \centering
  \begin{tabular}{cc}
    \includegraphics[width=60mm]{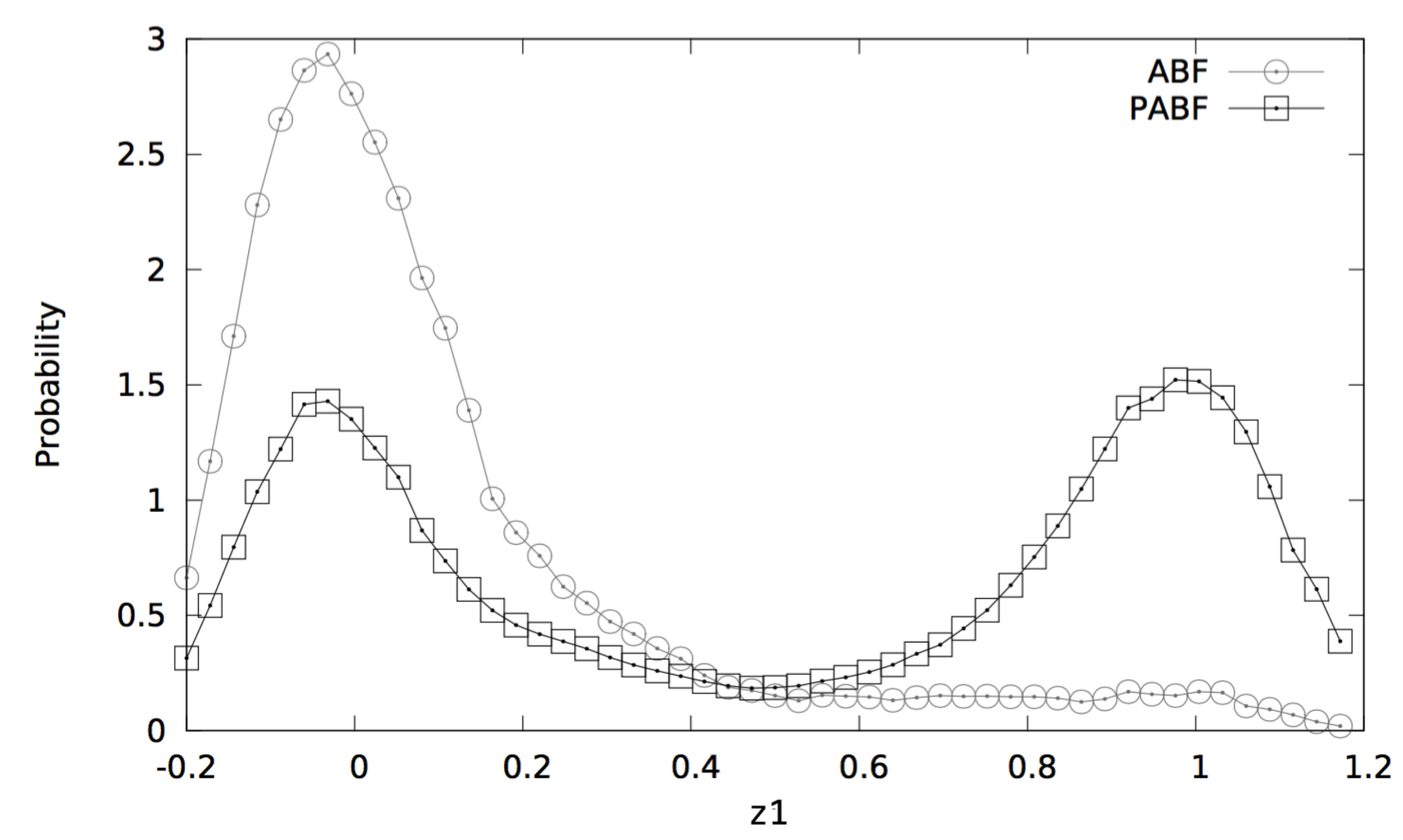}&
    \includegraphics[width=60mm]{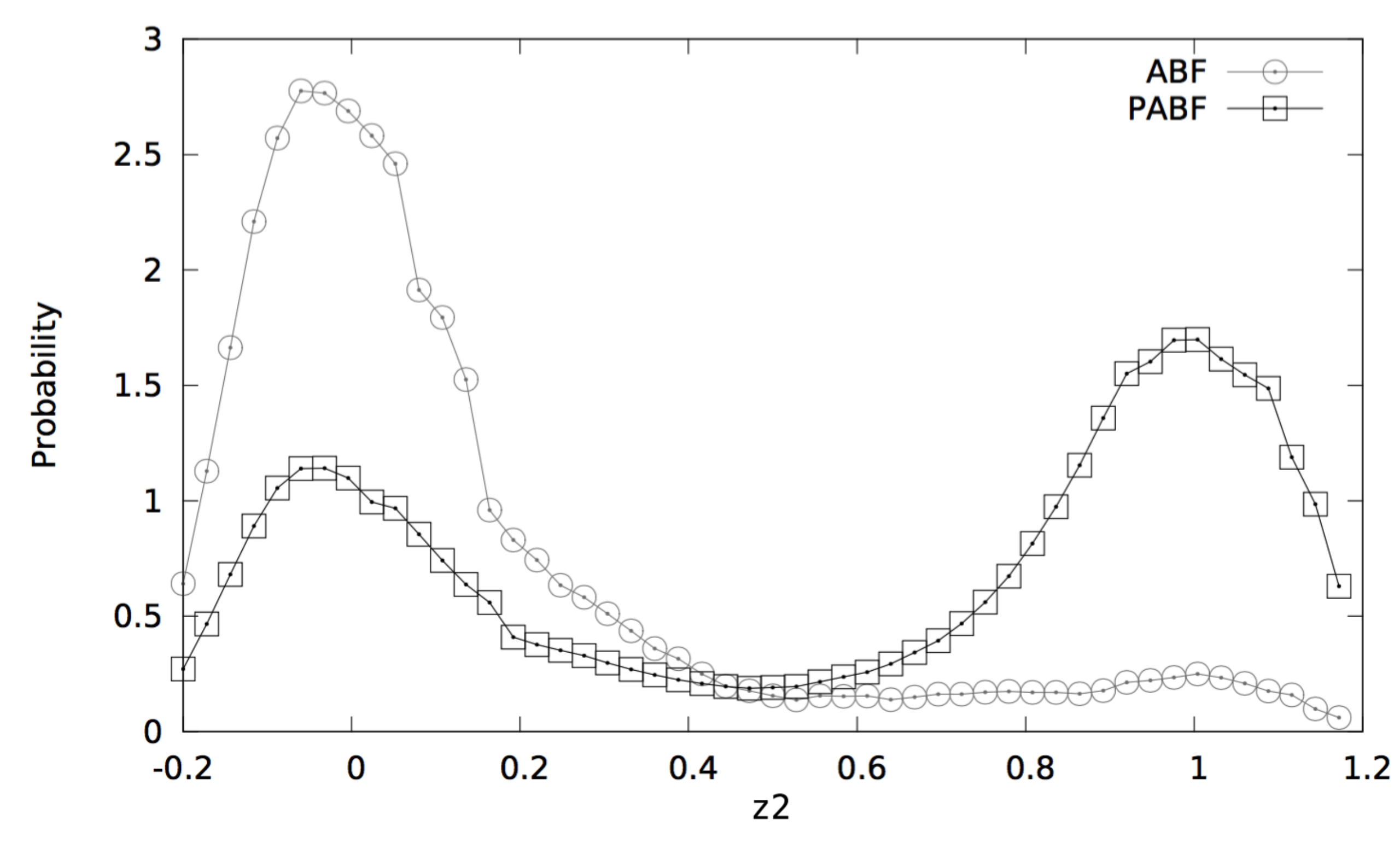}\\
  \end{tabular}
\caption{At time $5$. Left: $\int \psi^{\xi}(z_1,z_2)dz_2$; Right:  $\int \psi^{\xi}(z_1,z_2)dz_1$.}\label{dtr5}
\end{figure}

\begin{figure}[htp]
  \centering
  \begin{tabular}{cc}
    \includegraphics[width=60mm]{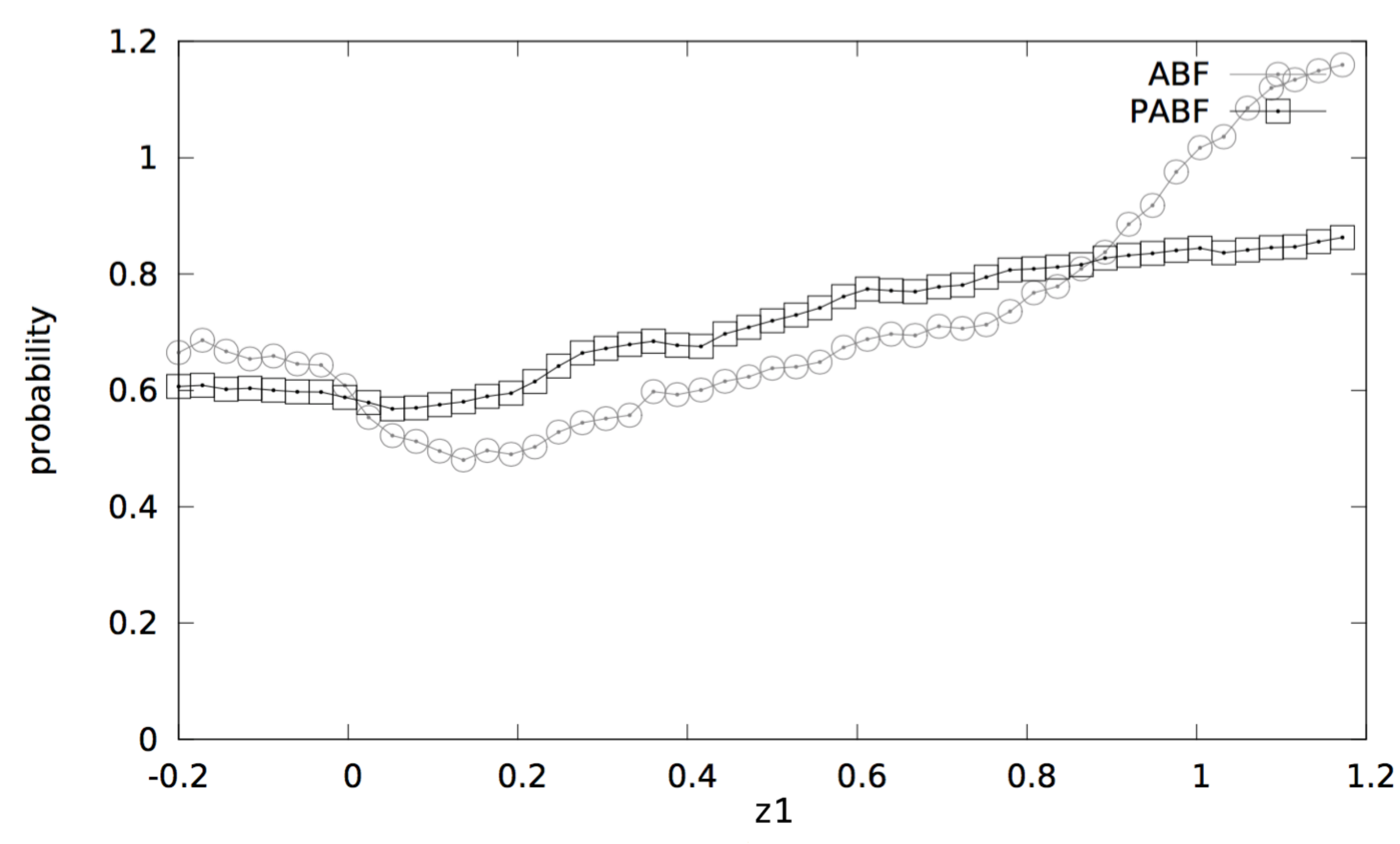}&
    \includegraphics[width=60mm]{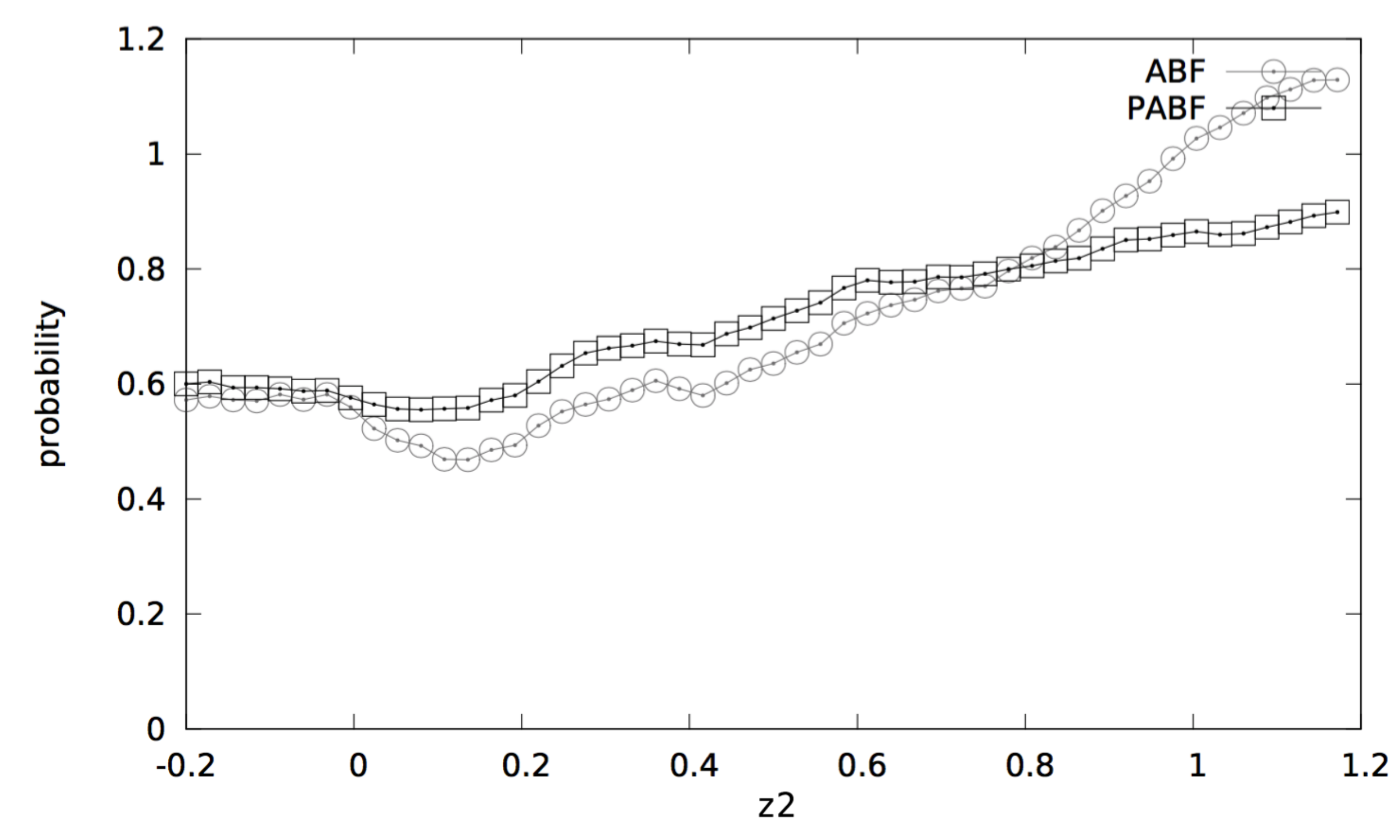}\\
  \end{tabular}
\caption{At time $25$. Left: $\int \psi^{\xi}(z_1,z_2)dz_2$; Right:  $\int \psi^{\xi}(z_1,z_2)dz_1$.}\label{dtr25}
\end{figure}

 Another way to illustrate that the projected ABF method converges faster than the standard ABF method is to plot the density function $\psi^{\xi}$ as a function of time. It is clearly observed (see Figures~\ref{dtr0}-\ref{dtr25}) that, for the PABF method, the convergence of $\psi^{\xi}$ to uniform law along $(\xi_1,\xi_2)$ is faster with the projected ABF method.

\newpage

%
%
%

\end{document}